\documentclass[runningheads]{llncs}
\usepackage{amsmath}
\usepackage{amssymb}
\usepackage{todonotes}
\usepackage{setspace}
\usepackage{overpic}
\usetikzlibrary{decorations.pathreplacing}

\DeclareMathOperator{\RR}{\mathbb{R}}
\newcommand{\sumx}{\sum_{i,j=0}^{N_x-1}}
\newcommand{\sums}{\sum_{p=0}^{N_s-1}}
\newcommand{\sumphi}{\sum_{q=0}^{N_\phi-1}}

\newcommand{\mres}{\mathbin{\vrule height 1.6ex depth 0pt width
0.13ex\vrule height 0.13ex depth 0pt width 1.3ex}}

\DeclareMathOperator{\NN}{\mathbb{N}}

\DeclareMathOperator{\Radon}{\mathcal{R}}
\DeclareMathOperator{\RayRadon}{\mathcal{R}^\mathrm{rd}_{\delta}}
\DeclareMathOperator{\PixelRadon}{\mathcal{R}^\mathrm{pd}_{\delta}}

\DeclareMathOperator{\RayRadonN}{\mathcal{R}^\mathrm{rd}_{\delta^n}}
\DeclareMathOperator{\PixelRadonN}{\mathcal{R}^\mathrm{pd}_{\delta^n}}

\DeclareMathOperator{\RayWeight}{\omega^\mathrm{rd}_{\delta}}
\DeclareMathOperator{\PixelWeight}{\omega^\mathrm{pd}_{\delta}}

\DeclareMathOperator{\imgdom}{\Omega}
\DeclareMathOperator{\sinodom}{\mathcal{S}}

\newcommand{\dd}[1]{\,\mathrm{d}{#1}}

\usepackage[T1]{fontenc}
\usepackage{graphicx}
\begin{document}
\title{A Novel Interpretation of the Radon Transform's Ray- and Pixel-Driven Discretizations under Balanced Resolutions}
\titlerunning{Novel interpretations of discretizations of the Radon transform}
\author{
Richard Huber\inst{1}\orcidID{0000-0003-1743-6786}
}
\authorrunning{R. Huber
}
\institute{Technical University of Denmark 
\\
Department of Applied Mathematics and Computer Science
\\
2800 Kgs. Lyngby, Denmark
\\
\email{richu@dtu.dk}
}
\maketitle              % typeset the header of the contribution
\begin{abstract}
Tomographic investigations are a central tool in medical applications, allowing doctors to image the interior of patients. The corresponding measurement process is commonly modeled by the Radon transform. In practice, the solution of the tomographic problem requires discretization of the Radon transform and its adjoint (called the backprojection). There are various discretization schemes; often structured around three discretization parameters: spatial-, detector-, and angular resolutions. The most widespread approach uses the ray-driven  Radon transform and the pixel-driven backprojection in a balanced resolution setting, i.e., the spatial resolution roughly equals the detector resolution. 
%Using these two different discretization frameworks leads to an unmatched operator pair, which can potentially harm the convergence of iterative solvers.
%Still, these methods are commonly used as they supposedly have better approximation properties, while their corresponding adjoints do not, i.e., the ray-driven forward is a much better approximation than the pixel-driven forward, and conversely for the backprojection is much more suitable than the ray-driven backprojection. 
The use of these particular discretization approaches is based on anecdotal reports of their approximation performance, but there is little rigorous analysis of these methods' approximation errors. This paper presents a novel interpretation of ray-driven and pixel-driven methods as convolutional discretizations, illustrating that from an abstract perspective these methods are similar. Moreover, we announce statements concerning the convergence  of the ray-driven Radon transform and the pixel-driven backprojection under balanced resolutions. Our considerations are supported by numerical experiments highlighting aspects of the discussed methods.

%This paper shows convergence in the strong operator topology (pointwise convergence) of the ray-driven Radon transform and the pixel-driven backprojection for balanced resolutions. 

\end{abstract}
\section{Introduction}

Computed Tomography (CT) is a crucial tool in medicine, allowing the investigation of the interior of patients' bodies \cite{Hsieh_CT_principles}.  The measurement process consists of emitting radiation that travels through the patient, the intensity of which is then measured by a detector on the other side.
This paper considers planar geometries (i.e., ignoring the third spatial dimension) with the patient's body located inside the two-dimensional unit ball ${\imgdom:=B(0,1)}\subset \RR^2$ (called the spatial domain).
The radiation is assumed to move along straight lines (beams)
\begin{equation}
L_{\phi,s}:=\{s \vartheta_\phi+t\vartheta_\phi^\perp \in \RR^2\big | \ t\in \RR\}
\end{equation}
 for  ${(\phi,s)\in\sinodom:=[0,\pi [ \times \mbox{]-1,1[}}$ (called the sinogram domain, having an angular and a detector component), where $\vartheta_\phi:=(\cos(\phi),\sin(\phi))\in \RR^2$ is the unit vector  associated with the angle $\phi$ and $\vartheta_\phi^\perp:=(-\sin(\phi),\cos(\phi))\in \RR^2$ is the direction rotated by 90 degrees counterclockwise; see Figure \ref{Fig_Geometry}.
The beam's loss of intensity is the accumulation of losses caused by the beam's interaction with the patient's mass within the body. 
This accumulation is commonly modeled by the  Radon transform \cite{Deans_Radon_applications_1993,Natterer:2001:MCT:500773} (in this context also called the forward projection), $\Radon \colon L^2(\imgdom)\to L^2(\sinodom)$ according to
\begin{align}
[\Radon f](\phi,s):= \int_{\imgdom} f(x) \dd {\mathcal{H}^{1}\mres L_{\phi,s}}(x) = \int_{-\sqrt{1-s^2}}^{\sqrt{1-s^2}} f(s\vartheta_\phi+t \vartheta_\phi^{\perp}) \dd t 
\end{align}
for $f\in L^2(\imgdom)$ and $(\phi,s)\in \sinodom$ (where $\mathcal{H}^1\mres L_{\phi,s}$ denotes the one-dimensional Hausdorff measure restricted to $L_{\phi,s}$), i.e., $\Radon$ is a collection of line integrals.
Thus, the tomographic reconstruction corresponds to the solution of the inverse problem $\Radon f =g$ for known measurements $g$ and unknown mass density distributions $f$.

While there is a direct inversion called the filtered backprojection \cite{Natterer:2001:MCT:500773}, many more evolved reconstruction techniques are iterative, e.g., iterative algebraic reconstruction algorithms \cite{SART_ALgo,Gilbert_Sirt,SCALES_CG_tomography_1987} and solution algorithms for convex optimization problems (e.g., total variation regularized reconstructions) \cite{Dong2013,C8NR09058K,Scherzer:2008:VMI:1502016}. These iterative methods also involve the adjoint operator $\Radon^*\colon L^2(\sinodom)\to L^2(\imgdom)$  (called the backprojection \cite{Natterer:2001:MCT:500773}), which,  given $g\in L^2(\sinodom)$, reads
\begin{equation} \label{equ_def_backprojection}
[\Radon^* g](x) := \int_{0}^{\pi} g(x\cdot \vartheta_\phi,\phi) \dd \phi  \qquad \text{for }x\in \imgdom.
\end{equation}

Only finite amounts of data can be measured and processed in practical applications. Thus, proper discretization $\Radon_\delta$ for $\Radon$ (and the adjoints) is imperative, where  $\delta$ describes the degree of discretization. For us, $\delta=(\delta_x,\delta_\phi,\delta_s)$ denotes the spatial resolution of reconstructions, and the angular and detector resolutions of measured data (we think of both reconstructions and measurements as images). The expectation is that with ever finer resolution $(\delta\to 0)$, also the approximation gets arbitrarily good (i.e., $\Radon_\delta \to \Radon$ in some sense).
The most widespread discretization framework uses the ray-driven Radon transform $\RayRadon$  \cite{doi:10.1118/1.4761867,Siddon1985FastCO,Path_through_pixels} and the pixel-driven backprojection $\PixelRadon^*$ \cite{doi:10.1137/20M1326635,4331812,Qiao2017ThreeNA,322963} (we call it an $\mathrm{rd}$-$\mathrm{pd}^*$ approach) with balanced resolutions, i.e., $\delta_x\approx \delta_s$. Note that usually $\delta_\phi$ and $\delta_s$ are determined from the measurements, while $\delta_x$ is somewhat selectable. The ray-driven approach discretizes line integrals by splitting them into sections associated with the intersections of the line with pixels, while the pixel-driven backprojection is based on a linear interpolation approximation of \eqref{equ_def_backprojection} (see Figure \ref{Fig_geometric_interpretation}).
There is also a ray-driven backprojection and a pixel-driven ray-driven Radon transform as the adjoint operators of the aforementioned operators. However, they do not see much use due to being supposed poor approximations that create artifacts \cite{7866886}.
Not using the forward and backprojections from the same framework (e.g., $\mathrm{rd}$-$\mathrm{rd}^*$ or $\mathrm{pd}$-$\mathrm{pd}^*$) can potentially harm iterative solvers \cite{10.1007/s00245-022-09933-5,870265}.

However, this danger seems to be outweighed in practice by the supposed better approximation performance of $\mathrm{rd}$-$\mathrm{pd}^*$ methods.
There is little rigorous analysis of approximation errors, and anecdotal knowledge of performance is much more prevalent.
\cite{doi:10.1137/20M1326635} rigorously discussed approximation errors for pixel-driven methods in the case that the spatial resolution $\delta_x$ is asymptotically smaller than the detector resolution $\delta_s$, finding convergence in the operator norm. Thus, a $\mathrm{pd}$-$\mathrm{pd}^*$ approach is justified when $\frac{\delta_x}{\delta_s}\to 0$ (both avoiding the unmatched operator issue, while yielding proper approximations). 
However, in practice, it is much more common to use balanced resolutions (i.e., $\delta_x\approx \delta_s$).
 This paper will partially justify the use of $\mathrm{rd}$-$\mathrm{pd}^*$ approaches for balanced resolutions  by 
 describing a novel interpretation of ray-driven methods that allows for much more structured analysis. In particular, based on these novel interpretations, we announce
 convergence results in the strong operator topology in Theorem \ref{Thm_approximation_ray_driven} (thus substantiating heuristic notions of approximation performance).
 The rigorous proof of the statement is available at arXiv \cite{Placeholder_archive}, with a corresponding journal article in preparation.
Some of these results were  already presented in the doctoral thesis \cite{huber2022pixel}.

\section{Discretization Approaches for the Radon transform}

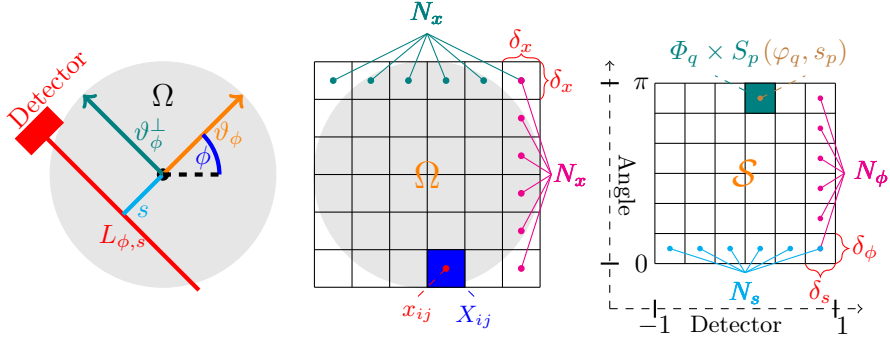
\begin{figure}
\centering
\usetikzlibrary{decorations.pathreplacing}
\usetikzlibrary{positioning,patterns}
\newcommand{\myfiguresize}{0.75}

\newcommand{\uppercutoff}{3.2}
\newcommand{\lowercutoff}{-3}
\newcommand{\mywidth}{ultra thick}
\begin{tikzpicture}[scale=\myfiguresize]
    \draw[gray!20,fill] (0,0) circle(2);
    \draw[fill] (0,0) circle (0.1);
    
    \draw (0,1.4) node{\large $\imgdom$};

    \clip(-2.8,\lowercutoff) rectangle (2,\uppercutoff);
    
    \pgfmathsetmacro{\myangle}{45}
    \pgfmathsetmacro{\beamend}{1.9}
    \pgfmathsetmacro{\detectorside}{0.4}
    \pgfmathsetmacro{\s}{1}
    \draw [red,rotate around ={\myangle+90:(0,0)}, \mywidth] (-\beamend,\s)--(\beamend,\s) node [midway, below] {$L_{\phi,s}$};
    
    \draw [cyan,rotate around ={\myangle+90:(0,0)}, \mywidth] (0,0)--(0,\s) node [midway, below] {$s$};
    
        \draw [fill,red ,rotate around ={\myangle+90:(0,0)}] (\beamend,\s-\detectorside)rectangle(\beamend+\detectorside,\s+\detectorside) node[above,yshift=0.5cm,xshift=0.6cm,rotate=45]{Detector};

    \draw [rotate around ={\myangle:(0,0)},orange,->, \mywidth] (0,0) -- (2,0) node[midway,right] {$\vartheta_\phi$};
    \draw [rotate around ={\myangle+90:(0,0)}, teal,->, \mywidth] (0,0) -- (2,0) node[midway,right] {$\vartheta_\phi^\perp$};

    \pgfmathsetmacro{\radius}{1}
    \draw [blue, \mywidth] (\radius,0) arc (0:\myangle:\radius) node [left,below] {$\phi$};
    \draw [dashed , \mywidth] (0,0) -- (\radius,0) ;
\end{tikzpicture}
\quad
\begin{tikzpicture}[scale=\myfiguresize]
  \draw[thick,fill, gray!20] (0,0) circle (2 );

    \clip(-2.05,\lowercutoff) rectangle (2.9,\uppercutoff);

\pgfmathsetmacro{\N}{6}
\foreach \n in {0,...,\N}
{
\pgfmathsetmacro{\s}{2*(\n/\N*2-1)}
\draw[] (-2,\s) -- (2,\s);
\draw[] (\s,-2) -- (\s,2);
}

\pgfmathsetmacro{\mydelta}{(1/\N*2)}
\foreach \n in {1,...,\N}
{
\pgfmathsetmacro{\s}{2*((\n-1/2)/\N*2-1)}
\draw[teal] (\s,2-\mydelta) -- (0,2.5) node [above] {$N_x$};
\draw[fill,teal] (\s,2-\mydelta) circle (0.05cm);

\draw[magenta] (2-\mydelta,\s) -- (2.2,0) node [right,xshift=-0.1cm] {$N_x$};
\draw[fill,magenta] (2-\mydelta,\s) circle (0.05cm);
}

\draw[red,decorate,decoration={brace,amplitude=3pt,mirror}] (2,2-2*\mydelta) -- (2,2) node[midway,right] {$\delta_x$};
\draw[red,decorate,decoration={brace,amplitude=3pt,mirror}] (2,2) -- (2-2*\mydelta,2) node[midway,above] {$\delta_x$};

\pgfmathsetmacro{\n}{4}
\pgfmathsetmacro{\m}{1}
\pgfmathsetmacro{\s}{2*((\n-1/2)/\N*2-1)}
\pgfmathsetmacro{\ss}{2*((\m-1/2)/\N*2-1)}

\draw[fill=blue] (\s-2*1/\N,\ss-2*1/\N) rectangle (\s+2*1/\N,\ss+2*1/\N);
\draw [dashed, blue] (\s,\ss) -- (\s+0.5,\ss-0.5)  node[below]{\footnotesize $X_{ij}$};

\draw [dashed, red] (\s,\ss) -- (\s-0.5,\ss-0.5)  node[below]{\footnotesize $x_{ij}$};
\draw [fill,red] (\s,\ss) circle (0.05);

\node[orange] at (0,0) {\Large $\imgdom$};

\end{tikzpicture}
\begin{tikzpicture}[scale=0.6]
\pgfmathsetmacro{\startpointx}{-3}
\pgfmathsetmacro{\startpointy}{-3}

\draw[dashed,->] (\startpointx,\startpointy) -- (\startpointx,2.5) node[midway,rotate=270,yshift=0.2cm]{\footnotesize Angle};
\draw[dashed,->] (\startpointx,-3) -- (2.6,-3)  node[midway,yshift=-0.2cm]{\footnotesize Detector};
\draw[thick] (\startpointx-0.2,2) --(\startpointx+0.2,2) node [xshift=0.3cm]{$\footnotesize \pi$};
\draw[thick] (\startpointx-0.2,-2) --(\startpointx+0.2,-2) node [xshift=0.3cm]{$\footnotesize 0$};
\draw[thick] (-2,-3+0.2) -- (-2,-3-0.2) node[yshift=-0.1cm]{$\footnotesize -1$};
\draw[thick] (2,-3+0.2) -- (2,-3-0.2) node[yshift=-0.1cm,xshift=0.1cm]{ $\footnotesize 1$};

\pgfmathsetmacro{\N}{6}
\foreach \n in {0,...,\N}
{
\pgfmathsetmacro{\s}{2*(2*\n/\N-1)}
\draw (-2,\s) -- (2,\s);
\draw (\s,-2) -- (\s,2);

}
\foreach \n in {1,...,\N}
{
\pgfmathsetmacro{\s}{2*(2*\n/\N-1)}
\draw[cyan](\s-2/\N, -2+2/\N)-- (0,-2.2) node[below]{$N_s$};
\draw[fill,cyan](\s-2/\N, -2+2/\N) circle (0.05);

\draw[magenta](2-2/\N,\s-2/\N)-- (2.2,0) node[right]{$N_\phi$};
\draw[fill,magenta](2-2/\N,\s-2/\N) circle (0.05);

}

\pgfmathsetmacro{\s}{2/\N}
\pgfmathsetmacro{\ss}{2-2/\N}

\draw[fill=teal](\s-2/\N,\ss-2/\N) rectangle (\s+2/\N,\ss+2/\N);
\draw[dashed,teal] (\s,\ss) -- (\s-1,\ss+0.5) node[above]{$\Phi_q \times S_p$};

\draw[dashed,brown] (\s,\ss) -- (\s+1,\ss+0.5) node[above]{$(\varphi_q,s_p)$};
\draw [fill,brown] (\s,\ss)  circle (0.05);

\draw [red,decorate,decoration={brace,amplitude=5pt,mirror}] (2,-2) -- (2,-2+4/\N) node[midway,right,xshift=.1cm] {$\delta_\phi$};

\draw [red, decorate,decoration={brace,amplitude=5pt,mirror}] (2-4/\N,-2) -- (2,-2) node[midway,below, yshift=-.1cm] {$\delta_s$};

\node[orange] at (0,0) {\Large $\sinodom$};

\end{tikzpicture}

    \caption{On the left is the geometry of the Radon transform; in the middle is the discretization of the spatial domain  $\imgdom$ into pixels $X_{ij}$ with width $\delta_x\times \delta_x$; on the right is the discretization of the sinogram domain $\sinodom$ into pixels $\Phi_q\times S_p$ of width $|\Phi_q|\times \delta_s$.}
    \label{Fig_Geometry}
\end{figure}

We start by discretizing the spatial domain $\imgdom=B(0,1)$ and the sinogram domain $\sinodom=[0,\pi[\times \mbox{]-1,1[}$ into `pixels'; one can think of data and reconstructions as digital two-dimensional images; see Figure \ref{Fig_Geometry}.
 We fix $N_x\in \NN$ and set $\delta_x=\frac{2}{N_x}$. We define the spatial pixel centers $x_{ij}:=( \frac{2i+1}{N_x}-1,\frac{2j+1}{N_x}-1)=(\left(i+\frac{1}{2}\right)\delta_x-1,\left(j+\frac{1}{2}\right)\delta_x-1)$ for $i,j\in [N_x]:= \{0,\dots,N_x-1\}$ and $X_{ij}:=x_{ij}+\left[-\frac {\delta_x}{2},\frac{\delta_x}{2}\right]^2$ denotes the corresponding squared spatial pixel with side-length $\delta_x$.
We consider a finite number $N_\phi\in \NN$ of angles $\phi_0< \dots <\phi_{N_\phi-1} \in \left[0,\pi \right[$ and associate them with the angular pixel $\Phi_q := \left[\frac{\phi_{q-1}+\phi_q}{2},\frac{\phi_{q+1}+\phi_q}{2}\right[$ for $q\in [N_\phi]$, where we understand $\phi_{-1}=\phi_{N_\phi-1}-\pi$ and $\phi_{N_\phi}=\phi_0+\pi$.
(Note that potentially $\Phi_0 \not \subset [0,\pi[$ or $\Phi_{N_\phi-1} \not \subset [0,\pi[$, in which case we tacitly $\pi$ periodically project them onto $[0,\pi[$.)
 Correspondingly, we set $\delta_\phi:=\max_{q\in [N_\phi]} |\Phi_q|$. For the sake of readability, we write $\vartheta_q$ when we mean $\vartheta_{\phi_q}$ for the unit vector associated with $\phi_q$. 
We assume a fixed number $N_s\in \NN$ of detector pixels and set $\delta_s=\frac{2}{N_s}$. 
The associated (equispaced) detector pixels are $S_p:= s_p+\left[-\frac{\delta_s}{2},\frac{\delta_s}{2}\right[$ for $p\in [N_s]$ with centers $s_p:= \frac{2p+1}{N_s}-1=(p+\frac 1 2 )\delta_s-1$. 
Hence, we have discretized the domain $\imgdom$ into a Cartesian $N_x\times N_x$ grid with resolution $\delta_x$, while the sinogram space $\sinodom$ contains $N_\phi\times N_s$ rectangular pixels $\Phi_q\times S_p$ for $q\in [N_\phi]$ and $p\in [N_s]$, i.e., with angular resolution $\delta_\phi$ and detector resolution $\delta_s$, see Figure \ref{Fig_Geometry}. We notationally combine all these resolutions to $\delta=(\delta_x,\delta_\phi,\delta_s)\in \RR^+\times\RR^+\times\RR^+$, and $N_x$, $N_\phi$ and $N_s$ are tacitly chosen accordingly.

One can naturally associate  pixels values of an image representing $f\in L^2(\imgdom)$ with average values $f_{ij}:=\frac{1}{\delta_x^2}\int_{X_{ij}}f(x) \dd x$. Moreover, one can consider the associated (pixelized) functions $f_\delta = \sumx f_{ij} u_{ij}$ in $U_\delta:=\text{span}\{u_{ij}\}\widehat = \RR^{N_x^2}$ with $u_{ij}:=\chi_{X_{ij}}-\frac 1 2 \chi_{\partial X_{ij}}$, where $\chi_M(x)$ equals 1 if $x\in M$  and zero otherwise, and  $\partial X_{ij}$ denotes the boundary of $X_{ij}$.
Analogously we define $V_\delta:=\text{span}\{\chi_{\Phi_q\times S_p}\}$ for discrete sinograms.
%Analogously, averages in sinogram pixels $g_{qp}:=\frac{1}{\delta_s|\Phi_q|}\int_{\Phi_q\times S_p}g(\phi,s) \dd{(\phi,s)}$ can be associated with functions $g_\delta\in V_\delta:=\text{span}\{\chi_{\Phi_q\times S_p}\}\widehat = \RR^{N_\phi N_s}$.
Discretizations of $\Radon$ translate to a matrix $A\in \RR^{ (N_\phi \cdot N_s) \times N_x^2 }$, where $A_{qpij}$ denotes the weight attributed to a pixel $X_{ij}$ in the calculation of $L_{\phi_q,s_p}$ ($q$ and $p$ determine a row in $A$, while $i$ and $j$ determine a column). Interpreting $A$ as a mapping from $U_\delta\subset L^2(\imgdom)$ to $V_\delta\subset L^2(\sinodom)$, we can understand the discretizations as finite rank operators between $L^2$ spaces, thus allowing for a proper comparison with $\Radon$. 

In order to compute $[\Radon f](\phi,s)$, one has to take relatively few values (the values of $f$ along $L_{\phi,s}$) into account; therefore also the $A_{qpij}$ should be non-zero only for pixels $X_{ij}$ that are close to $L_{\phi_q,s_p}$. To this end, we will set the weights $A_{qpij}$ using the following weight functions.   

\begin{definition}
Given $\delta$ and $\phi\in [0,\pi [$, we set $\overline s(\phi) := \frac{\delta_x}{2} (| \cos(\phi)|+| \sin(\phi)|)$, $\underline s(\phi) := \frac{\delta_x}{2} (\big | |\cos(\phi)|- |\sin(\phi)|\big|)$ and   $\kappa(\phi) :=  \min \left\{\frac{1}{|\cos(\phi)|},\frac{1}{|\sin(\phi)|}\right\}$ (with $\frac 1 0=\infty$). We define the ray-driven weight function for $t\in \RR$ according to
\begin{equation}\label{equ_def_ray_weight}
\RayWeight(\phi,t):= \frac{1}{\delta_x}
\begin{cases} 
 \frac{\overline s(\phi)-|t|}{\delta_x|\cos(\phi)\sin(\phi)|} \qquad & \text{if } |t| \in [\underline s(\phi),\overline s(\phi)[,
\\
\kappa(\phi) \qquad &\text{if } |t| \in <\underline s(\phi),
\\
\frac{1}{2} \qquad & \text{if } \phi\in \{ 0, \frac{\pi}{2}\}\text{ and } |t| =\overline s (\phi),
\\
0 & \text{else}.	
\end{cases}
\end{equation}
Moreover, we define the pixel-driven weight function to be
\begin{equation}
\PixelWeight(\phi,t)=\PixelWeight(t):= \frac{1}{\delta_s^2} \max\{\delta_s-|t|,0\} \qquad \text{for }\phi\in \left [0,\pi\right[, \  t\in \RR  .
\end{equation}
\end{definition}
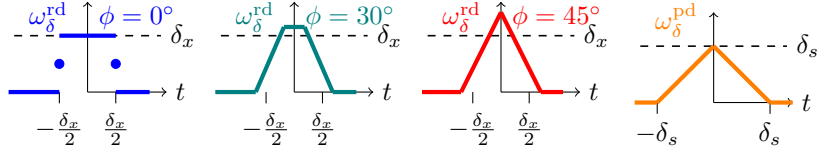
\begin{figure}
\centering
    \newcommand{\mywidthplot}{1.4}
    \newcommand{\myscale}{0.75}
    \newcommand{\myradius}{0.08}
    \newcommand{\mylinewidth}{ultra thick}
    \newcommand{\mycolor}{blue}
    \newcommand{\mytextpositionone}{-0.8}
    \newcommand{\mytextpositiontwo}{1}
    \newcommand{\myarrowup}{1.6}
    \begin{tikzpicture}

        \draw[scale=\myscale,->] (0,0) -- (\mywidthplot,0) node[right] {$t$};
        \draw[scale=\myscale,->] (0,0) -- (0,\myarrowup) ;

        \draw[scale=\myscale] (0.5,0) --(0.5,-0.2) node[below]{$\frac{\delta_x}{2}$};
        \draw[scale=\myscale] (-0.5,0) --(-0.5,-0.2) node[below]{$-\frac{\delta_x}{2}$};
        \draw[dashed, semithick, scale=\myscale] (-\mywidthplot+0.1,1) --(\mywidthplot-0.1,1) node[right]{$\delta_x$};

        \draw[scale=\myscale, domain=-\mywidthplot:-0.5, smooth, variable=\x, \mycolor,\mylinewidth] plot ({\x}, {0});
        \draw[scale=\myscale,fill,\mycolor] (-0.5,0.5) circle (\myradius cm);
        \draw[scale=\myscale, domain=-0.5:0.5, smooth, variable=\x, \mycolor,\mylinewidth] plot ({\x}, {1});
        \draw[scale=\myscale,fill,\mycolor] (0.5,0.5) circle (\myradius cm);
        \draw[scale=\myscale, domain=0.5:\mywidthplot-0.3, smooth, variable=\x, \mycolor,\mylinewidth] plot ({\x}, {0});

        \draw[\mycolor] (\mytextpositionone,\mytextpositiontwo) node[right,xshift=-0.1cm] {$\RayWeight \quad \phi=0^\circ$};;
        \end{tikzpicture}
\begin{tikzpicture}
    \newcommand{\mycolorthree}{teal}
\pgfmathsetmacro{\myangle}{60}
\pgfmathsetmacro{\supper}{(cos(\myangle)+sin(\myangle))/2}
\pgfmathsetmacro{\sunder}{(-cos(\myangle)+sin(\myangle))/2}
\pgfmathsetmacro{\mykappa}{1/sin(\myangle)}

        \draw[scale=\myscale,->] (0,0) -- (\mywidthplot,0) node[right] {$t$};
        \draw[scale=\myscale,->] (0,0) -- (0,\myarrowup) ;

        \draw[scale=\myscale] (0.5,0) --(0.5,-0.2) node[below]{$\frac{\delta_x}{2}$};
        \draw[scale=\myscale] (-0.5,0) --(-0.5,-0.2) node[below]{$-\frac{\delta_x}{2}$};
        \draw[dashed, semithick, scale=\myscale] (-\mywidthplot+0.1,1) --(\mywidthplot-0.1,1) node[right]{$\delta_x$};
    
\draw[scale=\myscale, domain=-\mywidthplot:-\supper, smooth, variable=\x, \mycolorthree,\mylinewidth] plot ({\x}, {0});
\draw[scale=\myscale, domain=-\sunder:\sunder, smooth, variable=\x, \mycolorthree,\mylinewidth] plot ({\x}, {\mykappa});
\draw[scale=\myscale, domain=-\supper:-\sunder, smooth, variable=\x, \mycolorthree,\mylinewidth] plot ({\x}, {\mykappa*(\x+\supper)/(\supper-\sunder)});
\draw[scale=\myscale, domain=\sunder:\supper, smooth, variable=\x, \mycolorthree,\mylinewidth] plot ({\x}, {\mykappa*(-\x+\supper)/(\supper-\sunder)});
\draw[scale=\myscale, domain=\supper:\mywidthplot-0.3, smooth, variable=\x, \mycolorthree,\mylinewidth] plot ({\x}, {0});

\draw[\mycolorthree] (\mytextpositionone,\mytextpositiontwo) node[right,xshift=-0.1cm] {$\RayWeight \quad \phi=30^\circ$};
\end{tikzpicture}
\begin{tikzpicture}
    \newcommand{\mycolorthree}{red}
    
\pgfmathsetmacro{\myangle}{45}
\pgfmathsetmacro{\supper}{(cos(\myangle)+sin(\myangle))/2}
\pgfmathsetmacro{\sunder}{(-cos(\myangle)+sin(\myangle))/2}
\pgfmathsetmacro{\mykappa}{1/sin(\myangle)}

        \draw[scale=\myscale,->] (0,0) -- (\mywidthplot,0) node[right] {$t$};
        \draw[scale=\myscale,->] (0,0) -- (0,\myarrowup) ;

        \draw[scale=\myscale] (0.5,0) --(0.5,-0.2) node[below]{$\frac{\delta_x}{2}$};
        \draw[scale=\myscale] (-0.5,0) --(-0.5,-0.2) node[below]{$-\frac{\delta_x}{2}$};
        \draw[dashed, semithick, scale=\myscale] (-\mywidthplot+0.1,1) --(\mywidthplot-0.1,1) node[right]{$\delta_x$};
    
\draw[scale=\myscale, domain=-\mywidthplot:-\supper, smooth, variable=\x, \mycolorthree,\mylinewidth] plot ({\x}, {0});
\draw[scale=\myscale, domain=-\sunder:\sunder, smooth, variable=\x, \mycolorthree,\mylinewidth] plot ({\x}, {\mykappa});
\draw[scale=\myscale, domain=-\supper:-\sunder, smooth, variable=\x, \mycolorthree,\mylinewidth] plot ({\x}, {\mykappa*(\x+\supper)/(\supper-\sunder)});
\draw[scale=\myscale, domain=\sunder:\supper, smooth, variable=\x, \mycolorthree,\mylinewidth] plot ({\x}, {\mykappa*(-\x+\supper)/(\supper-\sunder)});
\draw[scale=\myscale, domain=\supper:\mywidthplot-0.3, smooth, variable=\x, \mycolorthree,\mylinewidth] plot ({\x}, {0});

    \draw[\mycolorthree] (\mytextpositionone,\mytextpositiontwo) node[right,xshift=-0.1cm] {$\RayWeight \quad \phi=45^\circ$};
\end{tikzpicture}
\begin{tikzpicture}
    \newcommand{\mycolorfour}{orange}

        \draw[scale=\myscale,->] (0,0) -- (\mywidthplot,0) node[right] {$t$};
        \draw[scale=\myscale,->] (0,0) -- (0,\myarrowup) ;

        \draw[scale=\myscale] (1,0) --(1,-0.2) node[below]{$\delta_s$};
        \draw[scale=\myscale] (-1,0) --(-1,-0.2) node[below]{$-\delta_s$};
        \draw[dashed, semithick, scale=\myscale] (-\mywidthplot+0.1,1) --(\mywidthplot-0.1,1) node[right]{$\delta_s$};
    
\draw[scale=\myscale, domain=-\mywidthplot:-1, smooth, variable=\x, \mycolorfour,\mylinewidth] plot ({\x}, {0});
\draw[scale=\myscale, domain=-1:0, smooth, variable=\x, \mycolorfour,\mylinewidth] plot ({\x}, {1+\x});
\draw[scale=\myscale, domain=-0:1, smooth, variable=\x, \mycolorfour,\mylinewidth] plot ({\x}, {1-\x});

\draw[scale=\myscale, domain=1:\mywidthplot, smooth, variable=\x, \mycolorfour,\mylinewidth] plot ({\x}, {0});

    \draw[\mycolorfour] (\mytextpositionone,\mytextpositiontwo) node [right,xshift=-0.1cm,yshift=0.1cm] {$\PixelWeight$};
\end{tikzpicture}

\caption{Depiction of the ray-driven weight function $t\mapsto \delta_x^2\RayWeight(\phi,t)$ for fixed $\phi\in \{0^\circ,30^\circ,45^\circ\}$. On the right, the depiction of the pixel-driven weight $t\mapsto \delta_s^2\PixelWeight(t)$.}

\end{figure}

For the ray-driven method, one typically uses weights $A_{qpij} := \mathcal{H}^1(X_{ij}\cap L_{\phi_q,s_p})$ (i.e., the intersection length of line and pixel) computed in an iterative manner following the ray \cite{doi:10.1118/1.4761867}. The novel function $\delta_x^2 \RayWeight(\phi_q,x_{ij}\cdot\vartheta_q-s_p)$ is a closed form of this weight, allowing more structured analysis.

%%optional
% (The special case $\phi\in \{ 0,\frac{\pi}{2}\}\text{ and } |t| =\overline s (\phi)$ in \eqref{equ_def_ray_weight} relates to when a ray coincides with a side connecting two neighboring pixels, in which case we attribute half the intersection length to both pixels, but this choice is somewhat arbitrary.) 

%Both the ray-driven and pixel-driven methods can be understood as convolutional methods using the respective weight functions, only taking normal distances $x_{ij}\cdot\vartheta_q-s_p$ between $x_{ij}$ and $L_{\phi_q,s_p}$ into account, as we define next.

\begin{definition}
   Given a (weight) function $\omega \colon \sinodom \to \RR$, we define the convolutional Radon transform 
$\Radon_\omega\colon L^2(\imgdom)\to L^2(\sinodom)$, such that, for a function $f\in L^2(\imgdom)$, 
    \begin{equation}\label{equ_definition_convolutional_forward}
[\Radon_\omega f](\phi,s) := \sumphi\sums \chi_{\Phi_q \times S_p}(\phi,s) \sumx \omega (\phi_q,x_{ij}\cdot \vartheta_q -s_p) \int_{X_{ij}} f(x) \dd x.      
    \end{equation}
    The ray-driven Radon transform $\RayRadon$ and pixel-driven Radon transform $\PixelRadon$ are defined as special cases of $\Radon_\omega$ with $\omega=\RayWeight$ or $\omega=\PixelWeight$, respectively.
    The corresponding ray-driven or pixel-driven backprojections $\RayRadon^*$ and $\PixelRadon^*$ are  special cases of the convolutional backprojection $\Radon_\omega^*\colon L^2(\sinodom)\to L^2(\imgdom)$ according to
    \begin{equation}\label{equ_definition_convolutional_backprojection}
        [\Radon_\omega^* g](x) =   \sumx \chi_{X_{ij}}(x) \sumphi\sums \omega(\phi_q,x_{ij}\cdot \vartheta_q -s_p) \int_{\Phi_q\times S_p} g(\phi,s) \dd{(\phi,s)}
    \end{equation}
    for $g\in L^2(\sinodom)$ when setting $\omega$ to $\RayWeight$ or $\PixelWeight$, respectively.
\end{definition}

%Going to the discrete setting and understanding these discretizations as operators between $U_\delta$ and $V_\delta$ (or vice-versa) read as
%\begin{align}
%&[\RayRadon f]_{qp} = \delta_x^2\sumx \RayWeight(x_{ij}\cdot \vartheta_q -s_p,\phi_q)  f_{ij}
%\\
%&[\PixelRadon^* g]_{ij} = \sumphi\sums \delta_s |\Phi_q| \PixelWeight(x_{ij}\cdot \vartheta_q -s_p,\phi_q) g_{qp}.
%\end{align}

%The ray-driven forward and the pixel-driven backprojection have clear interpretations, as the weights attributed to a pixel by the ray-driven forward is the intersection length between pixels and lines, while the pixel-driven backprojection weights facilitate a linear interpolation; see Figure \ref{Fig_geometric_interpretation} and Lemma \ref{Lemma_weight_length}. 

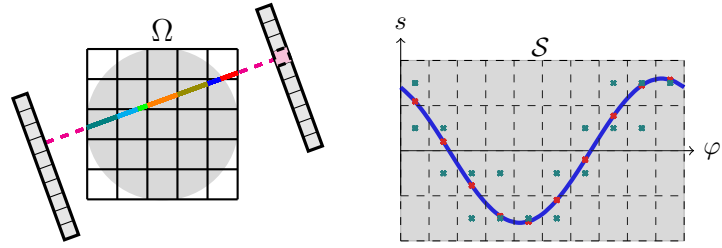
\begin{figure}
    \centering

\begin{tikzpicture}[scale=0.5]

 \definecolor{light_gray}{gray}{0.9}
\definecolor{light_gray2}{gray}{0.7}
\pgfmathsetmacro {\myangle}{20}
\pgfmathsetmacro {\myRE}{3.5}
\pgfmathsetmacro {\myRD}{3.4}
\pgfmathsetmacro {\myDW}{2}
\pgfmathsetmacro {\Ns}{10} 
\pgfmathsetmacro {\N}{5}

\pgfmathsetmacro {\myposition}{\Ns/2-2}
\pgfmathsetmacro {\myx}{\myRD+\myRE}
\pgfmathsetmacro {\myy}{-\myposition/\Ns*\myDW*2+\myDW-1/\Ns*\myDW}
\pgfmathsetmacro {\mypositionm}{\myposition-1}
\pgfmathsetmacro {\mypositionp}{\myposition+1}

\pgfmathsetmacro {\myyplus}{-\mypositionp/\Ns*\myDW*2+\myDW-1/\Ns*\myDW}
\pgfmathsetmacro {\myyminus}{-\mypositionm/\Ns*\myDW*2+\myDW-1/\Ns*\myDW}

\pgfmathsetmacro {\mycos}{cos(\myangle+90)};
\pgfmathsetmacro {\mysin}{sin(\myangle+90)};

\pgfmathsetmacro {\mynewx}{\mycos*\myx};
\pgfmathsetmacro {\mynewy}{\mysin*\myx};

\pgfmathsetmacro {\mycos}{cos(\myangle)};
\pgfmathsetmacro {\mysin}{sin(\myangle)};

\pgfmathsetmacro {\sourcex}{\mycos*(\myRE)-\mysin*(\myy)};
\pgfmathsetmacro {\sourcey}{\mysin*(\myRE)+ \mycos*\myy};

\pgfmathsetmacro {\mynorm}{sqrt(\mynewx*\mynewx+\mynewy*\mynewy)};
\pgfmathsetmacro {\mynewx}{\mynewx/\mynorm};
\pgfmathsetmacro {\mynewy}{\mynewy/\mynorm};

\pgfmathsetmacro {\projectionx}{\mynewy};
\pgfmathsetmacro {\projectiony}{-\mynewx};

\draw[] (0,2.5) node {\large $\imgdom$};
\draw[thick,gray!30,fill] (0,0) circle (2cm);

\foreach \x in {0,...,\N}
{
\draw[black,thick] (4/\N*\x-2,-2) -- (4/\N*\x-2,2);
\draw[black,thick] (-2,4/\N*\x-2) -- (2,4/\N*\x-2);
}

\foreach \x in {1,...,\N}
{
\foreach \y in {1,...,\N}
{
\pgfmathsetmacro {\posx}{4/\N*\x-2-2/\N};
\pgfmathsetmacro {\posy}{4/\N*\y-2-2/\N};

\pgfmathsetmacro {\posxnew}{\posx-\sourcex};
\pgfmathsetmacro {\posynew}{\posy-\sourcey};

\pgfmathsetmacro {\valp}{\posxnew*\mynewx+\posynew*\mynewy};

\draw[very thick,dashed, magenta,rotate around={\myangle:(0,0)}] (-\myRE,\myy) -- (\myRD+0.4,\myy) node [midway,cyan,xshift=-0.125cm,yshift=-0.2cm] {};

\newcommand{\mythickness}{ultra thick}

\pgfmathsetmacro {\val}{\posxnew*\mynewx+\posynew*\mynewy};
\pgfmathsetmacro {\valnormed}{\val};
\pgfmathsetmacro {\projx}{\posx-\valnormed*\mynewx};
\pgfmathsetmacro {\projy}{\posy-\valnormed*\mynewy};

\pgfmathsetmacro {\steigung}{\projectiony/\projectionx};
\pgfmathsetmacro {\startpointx}{2};
\pgfmathsetmacro {\startpointy}{\sourcey-\steigung*(\sourcex-\startpointx)};
%\draw[teal,\mythickness] (\sourcex,\sourcey) -- (\startpointx,\startpointy);

\pgfmathsetmacro {\t}{(\startpointy-6/5)/\steigung};
\pgfmathsetmacro {\newx}{\startpointx-\t};
\pgfmathsetmacro {\newy}{6/5};
\draw[red,\mythickness] (\newx,\newy) -- (\startpointx,\startpointy);
\pgfmathsetmacro {\startpointx}{\newx};
\pgfmathsetmacro {\startpointy}{\newy};

\pgfmathsetmacro {\t}{(\startpointx-6/5)};
\pgfmathsetmacro {\newx}{6/5};
\pgfmathsetmacro {\newy}{\startpointy-(\steigung*\t)};
\draw[blue,\mythickness] (\newx,\newy) -- (\startpointx,\startpointy);
\pgfmathsetmacro {\startpointx}{\newx};
\pgfmathsetmacro {\startpointy}{\newy};

\pgfmathsetmacro {\t}{(\startpointx-2/5)};
\pgfmathsetmacro {\newx}{2/5};
\pgfmathsetmacro {\newy}{\startpointy-(\steigung*\t)};
\draw[olive,\mythickness] (\newx,\newy) -- (\startpointx,\startpointy);
\pgfmathsetmacro {\startpointx}{\newx};
\pgfmathsetmacro {\startpointy}{\newy};

\pgfmathsetmacro {\t}{(\startpointx+2/5)};
\pgfmathsetmacro {\newx}{-2/5};
\pgfmathsetmacro {\newy}{\startpointy-(\steigung*\t)};
\draw[orange,\mythickness] (\newx,\newy) -- (\startpointx,\startpointy);
\pgfmathsetmacro {\startpointx}{\newx};
\pgfmathsetmacro {\startpointy}{\newy};

\pgfmathsetmacro {\t}{(\startpointy-2/5)/\steigung};
\pgfmathsetmacro {\newx}{\startpointx-\t};
\pgfmathsetmacro {\newy}{2/5};
\draw[green,\mythickness] (\newx,\newy) -- (\startpointx,\startpointy);
\pgfmathsetmacro {\startpointx}{\newx};
\pgfmathsetmacro {\startpointy}{\newy};

\pgfmathsetmacro {\t}{(\startpointx+6/5)};
\pgfmathsetmacro {\newx}{-6/5};
\pgfmathsetmacro {\newy}{\startpointy-(\steigung*\t)};
\draw[cyan,\mythickness] (\newx,\newy) -- (\startpointx,\startpointy);
\pgfmathsetmacro {\startpointx}{\newx};
\pgfmathsetmacro {\startpointy}{\newy};

\pgfmathsetmacro {\t}{(\startpointx+2)};
\pgfmathsetmacro {\newx}{-2};
\pgfmathsetmacro {\newy}{\startpointy-(\steigung*\t)};
\draw[teal,\mythickness] (\newx,\newy) -- (\startpointx,\startpointy);
\pgfmathsetmacro {\startpointx}{\newx};
\pgfmathsetmacro {\startpointy}{\newy};

}
}

\draw [very thick,fill=light_gray,rotate around={\myangle:(0,0)}] (\myRD,-\myDW) rectangle (\myRD+0.4,\myDW);

\draw [very thick,fill=light_gray,rotate around={\myangle:(0,0)}] (-\myRE,-\myDW) rectangle (-\myRE+0.4,\myDW);

\foreach \s in {1,...,\Ns}
{
\draw [rotate around={\myangle:(0,0)}] (\myRD,\s/\Ns*\myDW*2-\myDW) -- (\myRD+0.4,\s/\Ns*\myDW*2-\myDW) ;
\draw [rotate around={\myangle:(0,0)}] (-\myRE,\s/\Ns*\myDW*2-\myDW) -- (-\myRE+0.4,\s/\Ns*\myDW*2-\myDW);

}

%\draw[blue,thick,<-] %%decorate,decoration={brace,amplitude=5pt},xshift=0pt,yshift=0pt](1.69,1.69) -- (1.69-1.82,1.69+1.82) node [midway,xshift=-0.2cm,yshift=-0.2cm] {$\xi$};

\draw [very thick,dashed,fill=magenta!30!white,rotate around={\myangle:(0,0)}] (\myRD,-\myposition/\Ns*\myDW*2+\myDW-0/\Ns*\myDW*2) rectangle (\myRD+0.4,-\myposition/\Ns*\myDW*2+\myDW-1/\Ns*\myDW*2);
  \end{tikzpicture}
  \qquad
  \begin{tikzpicture}[scale=0.6]
\draw[->] (0, 0) -- (6.5, 0) node[right] {$\varphi$};
  \draw[->] (0, -1) -- (0, 2.5) node[above] {$s$};
   \pgfmathtruncatemacro\myindex{14}
  \pgfmathsetmacro {\mydomain}{6.28}
  \draw[scale=1, domain=0:\mydomain, smooth, variable=\x, blue, ultra thick] plot ({\x}, {2*0.8*cos(deg(\x+0.5))});
 
 \pgfmathtruncatemacro\Nangle{10}
  \foreach \i in {0,...,\Nangle}
  {
    \pgfmathsetmacro {\x}{6.28*\i/\Nangle}
  \draw[ dashed ] (\x,-2) -- (\x,2);
  }

 \pgfmathtruncatemacro\Ns{4}
  \foreach \i in {0,...,\Ns}
  {
    \pgfmathsetmacro {\x}{2*(2*\i/\Ns-1)}
  \draw[ dashed ] (0,\x) -- (6.28,\x);
  }

\foreach \i in {1,...,\Nangle} 
 {
     \pgfmathsetmacro {\x}{6.28*\i/\Nangle-6.28/\Nangle/2}
 \draw[ mark=x,scale=1, domain=\x:\x+0.01, smooth, variable=\x, red, thick] plot ({\x}, {2*0.8*cos(deg(\x+0.5))});
 }

 \begin{scope}

\clip (0,2) rectangle (6.28,-2);

 \foreach \i/\j in {1/3,2/2,3/1,4/1,5/0,6/1,7/2,8/3,9/3,10/4} 
 {
     \pgfmathsetmacro {\x}{6.28*\i/\Nangle-6.28/\Nangle/2}
   \pgfmathsetmacro {\y}{2*(2*\j/\Ns-1-1/\Ns)}
 \pgfmathsetmacro {\z}{2*(2*\j/\Ns-1+1/\Ns)}

\draw[ mark=x,scale=1, domain=\x:\x+0.01, smooth, variable=\x, teal, thick] plot ({\x}, {\y}); 
 
 \draw[ mark=x,scale=1, domain=\x:\x+0.01, smooth, variable=\x, teal, thick] plot ({\x}, {\z});
 }
 
  \end{scope}
  \draw [fill=gray,opacity=0.3] (0,2) rectangle (6.28,-2);
  \node[] at (3.1,2.3) {$\sinodom$};
\end{tikzpicture}
    \caption{Illustration of the ray-driven Radon transform (left) and pixel-driven backprojection (right). In the former, integration along a straight line is split into the sum of integrals along the intersections with pixels. The pixel-driven backprojection discretizes the angular integral \eqref{equ_def_backprojection} (along the blue curve $(x\cdot\vartheta_\phi,\phi)$) by a finite sum of angular evaluations $(x\cdot\vartheta_q,\phi_q)$ (red dots), with linear interpolation in the detector dimension (using the points represented by the green dots).}
    \label{Fig_geometric_interpretation}
\end{figure}

Note that these definitions are operators between $L^2$ spaces, whose corresponding matrices $A$ coincide with the classical interpretation of these methods in the literature as intersection lengths as weights and linear interpolation, respectively; see Figure \ref{Fig_geometric_interpretation}. Also, \eqref{equ_definition_convolutional_forward} and \eqref{equ_definition_convolutional_backprojection} show that from a certain point of view, the ray-driven and pixel-driven methods are actually quite similar. Therefore, they can also be implemented using structually analogous algorithms with similar degrees of parallelization.

\section{Theoretical Results}

\begin{lemma} \label{Lemma_weight_length}
Given $\delta$, and $f_\delta \in U_\delta$, for all $\hat q\in [N_\phi]$ and $\hat p\in [N_s]$, we have
\begin{equation}
\label{equ_lemma_exact_on_U}
[\Radon f_\delta] (\phi_{\hat q},s_{\hat p})= [\RayRadon f_\delta] (\phi_{\hat q},s_{\hat p}).  
\tag{$\mathrm{exact}^\mathrm{rd}$}
\end{equation}
Moreover, for fixed $\hat i,\hat j\in [N_x]$ such that $x_{\hat i \hat j}\cdot \vartheta_{\hat q}\in [s_0,s_{N_s-1}]$, we have
\begin{equation}
\label{equ_lemma_pixel_sum_s}
    \sums \PixelWeight(x_{\hat i \hat j}\cdot \vartheta_{\hat q} -s_p)= \frac{1}{\delta_s}
    \tag{$\mathrm{intpol}^{\mathrm {pd}}$}
\end{equation}
with exactly two non-zero summands $\hat p$ and $\hat p+1$ if $x_{\hat i \hat j}\cdot\vartheta_{\hat q}\in ]s_{\hat p},s_{\hat p+1}[$ and a single non-zero summand $\hat p$ if $x_{\hat i \hat j}\cdot\vartheta_{\hat q}=s_{\hat p}$.
\end{lemma}

One can interpret \eqref{equ_lemma_exact_on_U} as the ray-driven forward being precise for piecewise constant functions (images) in the detector pixel centers.
The ray-driven and pixel-driven discretizations have been known for decades, but no rigorous convergence results were presented. Theorem \ref{Thm_approximation_ray_driven} complements anecdotal reports on the performance of these discretization approaches, announcing convergence in the strong operator topology (the proof is found at \cite{Placeholder_archive}).

\begin{theorem} \label{Thm_approximation_ray_driven}
    Let $(\delta^n)_{n\in \NN}=(\delta_x^n,\delta_\phi^n,\delta_s^n)_{n\in \NN}$ be a sequence of discretization parameters with $\delta^n\overset{n\to \infty}{\to} 0$ (componentwise) and let $c>0$ be a constant. If $\frac{\delta_s^n}{\delta_x^n}<c$ for all $n\in \NN$, then,  for any $f\in L^2(\imgdom)$, we have
    \begin{equation}\label{equ_thm_ray_radon}
        \lim_{n\to \infty} \|\Radon f-\RayRadonN f\|_{L^2}=0.
        \tag{$\textrm{conv}^\mathrm{rd}$}
    \end{equation}
    If the sequence $(\delta^n)_{n\in \NN}$ satisfies $\frac{\delta_s^n}{\delta_x^n}\overset{n\to \infty}{\to} 0$, then, for each $g\in L^2(\sinodom)$, we have
        \begin{equation}\label{equ_thm_ray_estimate_backprojeciton}
        \lim_{n\to \infty} \|\Radon^* g-\RayRadonN^* g\|_{L^2}=0.
        \tag{$\textrm{conv}^\mathrm{rd*}$}
    \end{equation}
        If $\frac{\delta_x^n}{\delta_s^n} \leq c$ for all $n\in \NN$, then, for each $g\in L^2(\sinodom)$, we have
    \begin{equation}\label{equ_thm_pixel_backprojection}
        \lim_{n\to \infty}\|\Radon^* g - \PixelRadonN^* g\|_{L^2} =0
        \tag{$\textrm{conv}^\mathrm{pd*}$}.
    \end{equation}
    
\end{theorem}

\begin{remark}
Note that both \eqref{equ_thm_ray_radon} and \eqref{equ_thm_pixel_backprojection} are applicable in the case $\delta_x^n\approx \delta_s^n$ (or equivalently $N_x\approx N_s$), so using the $\mathrm{rd}$-$\mathrm{pd}^*$ approach in the balanced resolution case is justified in the sense that we have pointwise convergence of the operators.
Note that the convergence described in Theorem \ref{Thm_approximation_ray_driven} is not necessarily uniform, i.e., the convergence speed could depend on the specific $f$ and $g$.
\end{remark}

\section{Insights into Weight Functions via Simulations}

Throughout this section, we use small numerical experiments to highlight certain properties of the weight functions and show how they impact the discretized operators.  All calculations were implemented in Python using PyOpenCL with the pixel-driven  methods from Gratopy \cite{kristian_bredies_2021_5221443}, and with a custom ray-driven implementation in an analogous style. All calculations were executed on a 12th Gen Intel(R) Core(TM) i7-12650H in parallel with single precision. The code is available as a GitHub repository at \cite{Placeholder_GITHUB}.
We consider angles $\phi_q= \pi \frac{q}{N_\phi}$ for $q\in [N_\phi]$.
In our illustrations below, we always assume (and only visualize) $x,x_{ij}\in B(0,0.9)$ to avoid technical issues occuring at the boundary $\partial \imgdom$ that are not central to our observations.

\subsection*{Example 1}
Let us start with the most trivial example for backprojections; we consider the constant sinogram  $g(\phi,s)=1$ for $(\phi,s)\in \sinodom$.
It is trivial to see that the backprojection satisfies $[\Radon^* g](x)=\pi$ (as the integrand in \eqref{equ_def_backprojection} is constant). When looking at \eqref{equ_definition_convolutional_backprojection} for this concrete example, we would like to have 
\begin{equation} \label{equ_summability_all_angles}
   \sumphi \sums \delta_s |\phi_q|\omega(\phi_q,x_{ij}\cdot \vartheta_q -s_p) \underset{\text{def}}{\overset{\text{per }}{=}} [\Radon_\omega^* g](x_{ij})=[\Radon^*g](x_{ij})=\pi.
\end{equation}
(Requesting \eqref{equ_summability_all_angles} to hold for all $x\in \imgdom$ seems quite naturally then.)
When \eqref{equ_summability_all_angles} is satisfied, the method will be precise for this concrete example. And the pixel-driven backprojection (i.e., with $\omega=\PixelWeight$) satisfies \eqref{equ_summability_all_angles}  due to \eqref{equ_lemma_pixel_sum_s} (for $\|x\|\leq |s_0|$) and since $\sumphi |\Phi_q|=\pi$. Hence, the pixel-driven backprojection should be correct up to machine precision ($\approx 10^{-6}$) for this concrete example, which numerical experiments confirm. 

\begin{figure}
\newcommand{\mycolor}{cyan}
\newcommand{\mywidth}{0.22}
    \centering
   \textbf{Approximation Errors of Ray-Driven Backprojections for Example 1}\\
\begin{tabular}{ccc}
 \begin{overpic}[height=\mywidth\textwidth]{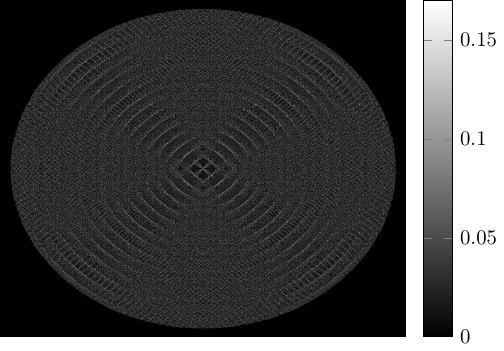}
\put(05,65){\color{\mycolor} \tiny a) 
\begin{minipage}[t]{0.3\textwidth}
$N_s=1000$, $N_x=1000$, \\
$N_\phi=90$\end{minipage} }
\put(20,05){\color{\mycolor} \tiny  rel. Error $0.012$}
\end{overpic}    

&

    \begin{overpic}[height=\mywidth\textwidth]{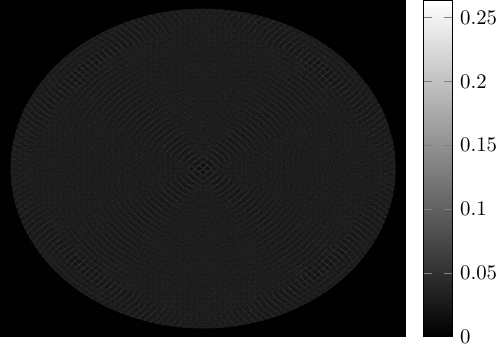}
\put(05,65){\color{\mycolor} \tiny b) 
\begin{minipage}[t]{0.3\textwidth}
$N_s=2000$, $N_x=2000$, \\
$N_\phi=90$\end{minipage} }
\put(20,05){\color{\mycolor} \tiny  rel. Error $0.012$}
    \end{overpic}

    &

    \begin{overpic}[height=\mywidth\textwidth]{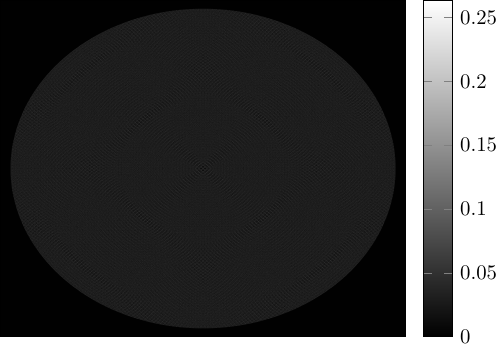}
\put(05,65){\color{\mycolor} \tiny c) 
\begin{minipage}[t]{0.3\textwidth}
$N_s=4000$, $N_x=4000$, \\
$N_\phi=90$\end{minipage} }
\put(20,05){\color{\mycolor} \tiny  rel. Error $0.012$}
    \end{overpic}

\\    
    
    \begin{overpic}[height=\mywidth\textwidth]{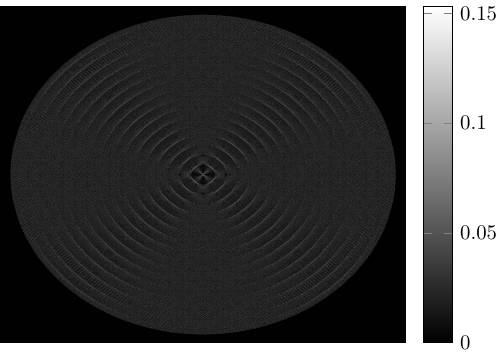}
\put(05,65){\color{\mycolor} \tiny d) 
\begin{minipage}[t]{0.3\textwidth}
$N_s=2000$, $N_x=2000$, \\
$N_\phi=180$\end{minipage} }
\put(20,05){\color{\mycolor} \tiny  rel. Error $0.0086$}
\end{overpic}    

&

    \begin{overpic}[height=\mywidth\textwidth]{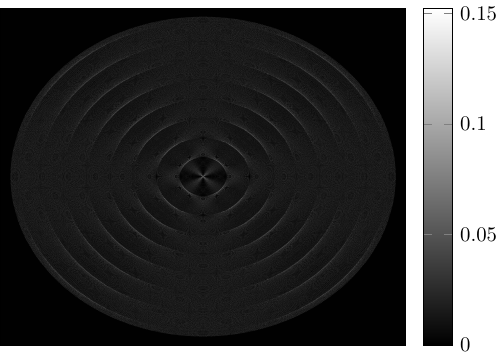}
\put(05,65){\color{\mycolor} \tiny e) 
\begin{minipage}[t]{0.3\textwidth}
$N_s=2000$, $N_x=2000$, \\
$N_\phi=360$\end{minipage} }
\put(20,05){\color{\mycolor} \tiny  rel. Error $0.0061$}
    \end{overpic}
    
    &
    
    \begin{overpic}[height=\mywidth\textwidth]{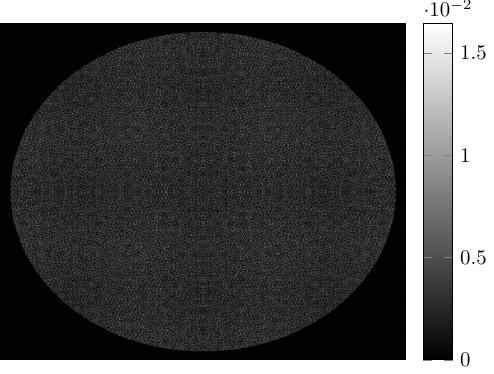}
\put(05,65){\color{\mycolor} \tiny f) 
\begin{minipage}[t]{0.3\textwidth}
$N_s=4000$, $N_x=1000$, \\
$N_\phi=90$\end{minipage} }
\put(20,05){\color{\mycolor} \tiny  rel. Error $0.0011$}
    \end{overpic}
\end{tabular}      
\caption{Illustration of the error incurred by ray-driven backprojections for Example 1. In the first row, a balanced resolutions setting with fixed $N_\phi=90$ but increasing $N_x=N_s$ is shown. The amplitude of the errors and the relative $L^2$ error, do not reduce with finer spatial and detector resolutions. In d) and e), we increase $N_\phi$ to $N_\phi=180$ and $N_\phi=360$, yielding slight improvements. Lastly, we depict the error when reducing the spatial resolution to $N_x=1000$ while keeping $N_s=4000$ detector pixels, leaving the balanced resolution setting and significantly reducing the error.}
 
    \label{Fig_example_1}

\end{figure}
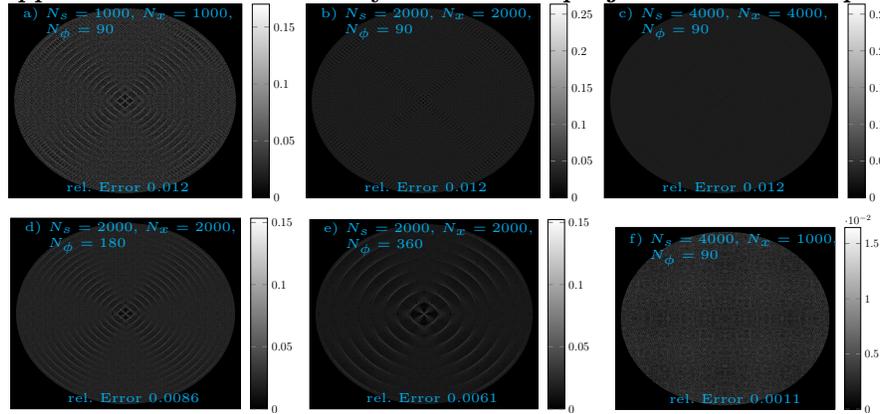

In contrast, for the ray-driven methods, it is unclear whether \eqref{equ_summability_all_angles} holds. 
Figures \ref{Fig_example_1} a), b) and c) depict the error of simulations using the ray-driven backprojections, revealing some geometric patterns whose amplitude does not decrease with increased $N_x$ and $N_s$ under balanced resolutions (i.e., $N_x=N_s$). Also, the corresponding $L^2$ error stays virtually the same.  
The relative error is roughly one percent, i.e., not huge. Still, this illustrates a potential issue that appears to not be solved by increasing resolutions. Note that the structured nature of the error suggests this is not caused by random rounding errors, but has a deeper cause in the weight function $\RayWeight$. 

Looking at Figures \ref{Fig_example_1} c), d) and e), we note that the error appears to reduce with increasing number of angles, though not very quickly. 
As illustrated in the pixel-driven backprojections' precision,  the angular resolution should not be a natural limiting factor for this example. In spite of that, for the ray-driven method it appears so. This can be understood as follows, \eqref{equ_summability_all_angles} is a mean value of $\phi \mapsto \sums  \delta_s \omega(\phi,x_{ij}\cdot \vartheta_\phi -s_p)$, and while this does not have a fixed value for all angles $\phi$ or points $x_{ij}$, it still averages out to $1$. Thus, the more angles, are considered the closer is \eqref{equ_summability_all_angles} to be satisfied.

According to \eqref{equ_thm_ray_estimate_backprojeciton}, if $\delta_x\ll \delta_s$ the method should converge. And indeed, as depicted in Figure \ref{Fig_example_1} f), the error reduces significantly (by a factor of 10) when considering $N_x=1000$ and $N_s=4000$. It might seem counter-intuitive that using coarser spatial resolution (compared to \ref{Fig_example_1} c)) yields a significantly better approximation, but this is in line with \eqref{equ_thm_ray_estimate_backprojeciton}.

%%optional
%So even though this is (in some sense) the easiest imaginable example (a constant sinogram), this already gave us quite a bit of insight into requirements of weight functions via \eqref{equ_summability_all_angles}.

\subsection*{Example 2}

 We fix some $\hat q\in [N_\phi]$ and consider $g_{\hat q}(\phi,s)=\frac{1}{|\Phi_{\hat q}|}$ for ${ s \in ]-1,1[}$ and $\phi \in \Phi_{\hat q}$   and zero otherwise, i.e., the sinogram has only very small angular support.  Straight-forward calculation shows $[\Radon^* g_{\hat q}](x) =1$ for all $x\in \imgdom$. 

Hence, in order for the convolutional approximations of this example (with the specific $\hat q$) to be precise, one would require 
\begin{equation}\label{equ_summability_single_angle}
\sums  \delta_s \omega(x_{ij}\cdot \vartheta_{\hat q} -s_p,\phi_{\hat q}) =1.
\end{equation}
(Naturally, $\hat q$ was arbitrary, so we would want \eqref{equ_summability_single_angle} to hold for all angles.)
 For the pixel-driven method, \eqref{equ_summability_single_angle} holds via \eqref{equ_lemma_pixel_sum_s}, and therefore, the pixel-driven method is again precise (up to machine precision).

The averaging with respect to the angles described in Example 1 that moderated \eqref{equ_summability_all_angles} does not apply here.
Therefore, the issue seen in Example 1 becomes much more severe for the ray-driven method.
In fact, a key step in proving  convergence of the ray-driven backprojection in \eqref{equ_thm_ray_estimate_backprojeciton} is to show that ray-driven weight functions asymptotically satisfy \eqref{equ_summability_single_angle} when $\delta_s \ll \delta_x$.

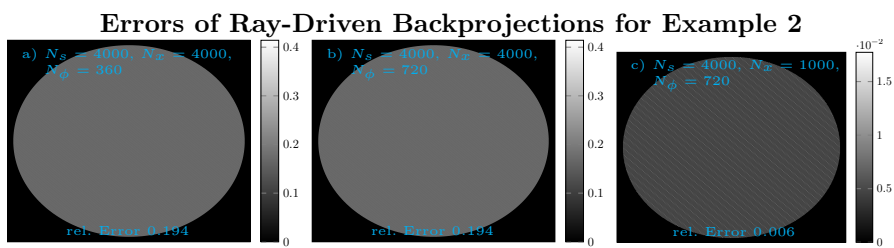
\begin{figure}
\centering
\textbf{Errors of Ray-Driven Backprojections for Example 2}
\newcommand{\mycolor}{cyan}
\newcommand{\mywidth}{0.23}
\begin{overpic}[height=\mywidth\textwidth]{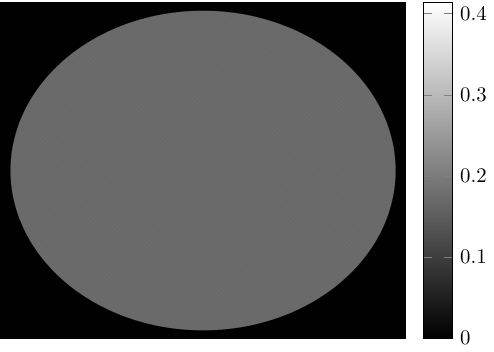}
\put(05,65){\color{\mycolor} \tiny a) 
\begin{minipage}[t]{0.3\textwidth}
$N_s=4000$, $N_x=4000$, \\
$N_\phi=360$\end{minipage} }
\put(20,05){\color{\mycolor} \tiny  rel. Error $0.194$}
\end{overpic}
\begin{overpic}[height=\mywidth\textwidth]{Backprojection/piecewise_constant/720/Differences_ray_Ns4000_Nx4000_error=0.1938212524397751.pdf}
\put(05,65){\color{\mycolor} \tiny b) 
\begin{minipage}[t]{0.3\textwidth}
$N_s=4000$, $N_x=4000$, \\
$N_\phi=720$\end{minipage} }
\put(20,05){\color{\mycolor} \tiny  rel. Error $0.194$}
\end{overpic}
\begin{overpic}[height=\mywidth\textwidth]{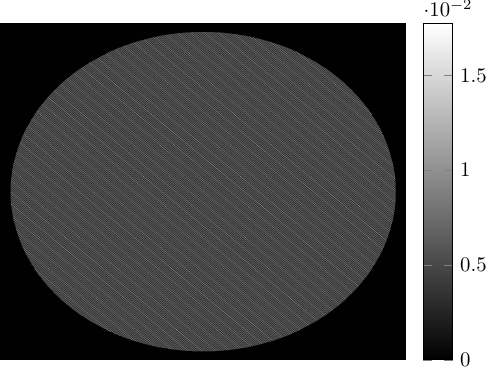}
\put(05,65){\color{\mycolor} \tiny c) 
\begin{minipage}[t]{0.3\textwidth}
$N_s=4000$, $N_x=1000$, \\
$N_\phi=720$\end{minipage} }
\put(20,05){\color{\mycolor} \tiny  rel. Error $0.006$}
\end{overpic}
\caption{Illustration of errors for Example 2, depicting the approximation errors of ray-driven backprojection with $\hat q$ such that $\phi_{\hat q}=\frac \pi 4$. Under balanced resolutions in a) and b), these errors are very significant and do not reduce with an increased number of angles. In c), we leave the balanced resolution setting, reducing errors significantly.}
\label{Fig_example2}

\end{figure}

Figures \ref{Fig_example2} a) and b) depict the errors of the ray-driven backprojection (with $\hat q$ such that $\phi_{\hat q}=\frac{\pi}{4}$), showing again certain patterns
which are much more prominent than they were in Example 1. They cause a relative $L^2$ error of 19\% and again appear to remain while in the balanced resolutions setting with finer resolutions.
However, increasing the number of angles does not improve things this time.
%%optional
% The error appears uniform everywhere, but appon\todo{uppon?} closer inspection it has a very fine structure. 
Again, choosing $N_x$ significantly smaller than $N_s$ reduces this error significantly (by a factor of roughly 30); see Figure \ref{Fig_example2} c).

It might be possible to translate these considerations into a rigorous construction of a function $g$ (only being non-zero on tacitly chosen angles) with $\| \Radon^*g - \RayRadon^* g\| \geq c \|g\|$ for some constant $c>0$ and all $N_x,N_s$ such that $N_x \approx N_s$ and $\delta_\phi \to 0$. If that was  the case, the ray-driven Radon transform could not converge to the Radon transform in the operator norm. Then, for any $\delta$ there would be a function $f=f(\delta)$ such that $\|\Radon f -\RayRadon f\|\geq c \|f\|$, i.e., no uniform convergence would hold and no matter how fine the resolution, there would always be functions for which the ray-driven Radon transform is a poor approximation.

\subsection*{Example 3}
Next, we consider the sinogram $g(\phi,s)= s$ for all $\phi\in[0,\pi[$ and $s\in \RR$ with 
\begin{equation}
\label{equ_example_3_integral}
[\Radon^* g](x) = \int_0^\pi \bold x \cos(\phi)+ \bold y \sin(\phi) \dd \phi= 2 \bold y ,
\end{equation}
 where $x=(\bold x, \bold y)$. Since the pixel-driven method is a linear interpolation (see \eqref{equ_lemma_pixel_sum_s}), and $g$ is linear in the detector dimension and constant in the angular dimension, one might suspect that $\PixelRadon^*$ is again accurate up to machine precision. However, this would be a fallacy. One can explicitly calculate the pixel-driven backprojection to be
\begin{equation}
[\PixelRadon^* g](x_{ij}) = \sumphi \bold x_i \cos(\phi_q) + \bold y_j \sin(\phi_q),
\end{equation}
where $x_{ij}=(\bold x_i,\bold y_j)$.
Hence, the error that the pixel-driven backprojection suffers is the error of the Riemann sum of the integrals in \eqref{equ_example_3_integral}, i.e., 
\begin{equation}
[\Radon^*g](x_{ij})-[\PixelRadon g](x_{ij})= -\bold x_i \sumphi \cos(\phi_q)+ \bold y_j \left(2-\sumphi\sin(\phi_q)\right).
\end{equation}
 
Thus, the pixel-driven backprojection can lead to highly structured errors related to how well these Riemann sums approximate the cosine and sine integrals. In particular, although $g$ does not depend on $\phi$, the angular resolution can have an impact here.

\begin{figure}
\newcommand{\mywidth}{0.25}
\newcommand{\mycolor}{green}
\centering
\textbf{Backprojection Errors in Example 3}\\
\begin{overpic}[height=\mywidth\textwidth]{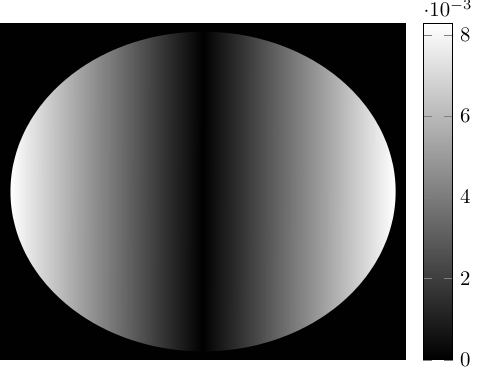}
\put(05,65){\color{\mycolor} \tiny a) 
\begin{minipage}[t]{0.3\textwidth}Naive angles\\
$N_s=4000$, $N_x=4000$, \\
$N_\phi=360$\end{minipage} }
\put(20,05){\color{\mycolor} \tiny  rel. Error $0.0043$}
\end{overpic}
\begin{overpic}[height=\mywidth\textwidth]{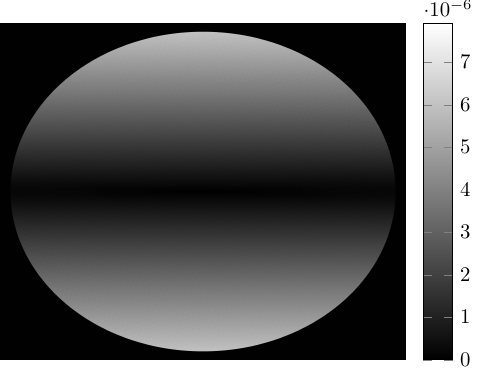}
\put(05,65){\color{\mycolor} \tiny b) 
\begin{minipage}[t]{0.3\textwidth}Shifted angles\\
$N_s=4000$, $N_x=4000$, \\
$N_\phi=360$\end{minipage} }
\put(20,05){\color{\mycolor} \tiny  rel. Error $3.18 \ 10^{-6}$}
\end{overpic}
\begin{overpic}[height=\mywidth\textwidth]{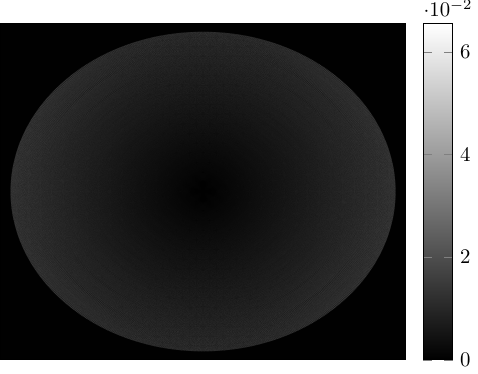}
\put(05,65){\color{\mycolor} \tiny c) 
\begin{minipage}[t]{0.3\textwidth}
$N_s=4000$, $N_x=4000$, \\
$N_\phi=360$\end{minipage} }
\put(20,05){\color{\mycolor} \tiny  rel. Error $0.0106$}
\end{overpic}

\caption{Illustration of the error of backprojections for Example 3.  In a) the pixel-driven backprojections' errors with the naive angles, and in b) with the described angle shift. The ray-driven backprojection's error is depicted in c). As can be seen the pixel-driven methods error significantly depends on the angles. %For the ray-driven backprojection, the improvement from the angular shift is much lesser, so other factors suposedly drive the error.
}
\label{Fig_example_3}
\end{figure}
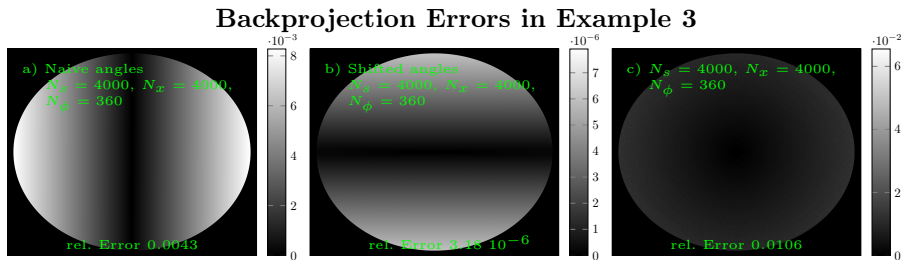

For our numerical experiments, we started with the naive choice $\phi_0=0^\circ$, $\phi_1=1^\circ$, \dots, $\phi_{179}=179^\circ$. Then, the error of the Riemann sum approximating the cosine integral is roughly $\pi/180 \approx 0.0174$ while the Riemann sum of the sine terms equals roughly $1.99995$. Hence, when computing the pixel-driven method for this example, the error is dominated by the linear $0.0174\ \bold x_i$, see Figure \ref{Fig_example_3} a).

 If we shift the considered angles by $0.5^\circ$, i.e., $\phi_0=0.5^\circ$, $\phi_1=1.5^\circ$, and so on, by symmetry the cosine Rieman sum is zero, while the sine Riemann sum's error is roughly $2\ 10^{-5}$. 
 Hence, the error is dominated by $2 \bold y_j 10^{-5}$, which overall reduces the approximation error substantially (roughly by a factor of $10^{-3}$). This structure is also observed in our simulations in Figure \ref{Fig_example_3} b). Note that neighter of the presented backprojections is terrible, however, their difference and the structure of the errors is unexpected.

Naturally, for real data, such considerations cannot be taken into account (as we have no insight into the functions we are integrating). Nonetheless, this gives some interesting insight that even though sinograms are not changing with respect to the angular domain, the angular discretization can significantly impact approximation errors.
Theorem \ref{Thm_approximation_ray_driven} suggests that with increasing number of angles this issue reduces (irrespective of the specific angles).
In contrast, in the ray-driven backprojection the error is not as structured (see Figure \ref{Fig_example_3} c)), suggesting that factors other than the Riemann sum drive the error.

%\subsection*{Example 4}
%
%We consider the function $g(\phi,s)=s^2*\sin(\phi)$ for $s\in ]-1,1[$ and $\phi \in [-\frac{\pi}{2},\frac{\pi}{2}[$. Straight forward calculation shows that $[\Radon^* g](x)==x_1x_2\frac{2}{3}$. Note that this is now a quadratic function, thus here also the pixel-driven method will incurr an interpolation error in the 
%
%One can observe similar issues with the angular resolution creating certain errors. To reduce them, we consider a very high number of angles $N_\phi=10000$. As can be seen

\section{Conclusion and Outlook}
This paper presented a novel interpretation of ray-driven and pixel-driven methods as `convolutional discretizations', suggesting in some sense the methods are more alike than previously thought.
Moreover, we announced novel convergence statements in the strong operator topology in Theorem \ref{Thm_approximation_ray_driven}. These results give a theoretical foundation for the widespread use of ray-driven forward and pixel-driven backprojection operators under balanced resolutions (previously based on anecdotal approximation properties).  These anecdotes also suggest that ray-driven backprojections and pixel-driven forward projections do not converge. While certainly no proof, our simulations support the former. Hence, the presented results might suggest that the combination of ray-driven and pixel-driven methods is the best one can do while maintaining the balanced resolution setting. 

The convergence results are probably extendable to other types of tomography (with fanbeam or conebeam operators), and might be the topic of future investigations.
Note that this work did not address the issue of unmatched operators. However, unmatched operators may be preferable if the alternative uses non-approximating discretizations.   

The convergence result \eqref{equ_thm_ray_estimate_backprojeciton} shows that while the ray-driven backprojection might be unsuitable for balanced resolutions, it approximates the backprojection if $\frac{\delta_s}{\delta_x}\to 0$. Thus, future works might investigate if convergence in the operator norm is achieved in that setting. Moreover, Example 2 might yield a way to show that ray-driven methods do not converge in the operator norm when using balanced resolutions, and thus might also merit future investigations.

\begin{credits}
\subsubsection{\ackname}
RH was partially supported by The Villum Foundation (Grant No.25893) and partially by the International Research Training Group ``Optimization and Numerical Analysis for Partial Differential Equations with Nonsmooth Structures,"" funded by the German Research Council (DFG) and Austrian Science Fund (FWF) grant W1244.

%\subsubsection{\discintname}
%The authors have no competing interests to declare that are relevant to the content of this article. 
\end{credits}
 \bibliographystyle{SSVM_minipaper}
\bibliography{splncs04}

\begin{thebibliography}{10}
\bibitem{Placeholder_GITHUB}
*: This is a placeholder for the github repository not yet uploaded containing
  the code with which the simulations where performed. GITHUB  (2025)
\end{thebibliography}

\end{document}